\documentclass[smallcondensed]{svjour3}
\usepackage[leqno]{amsmath}
\usepackage[all]{xy}
\usepackage[mathscr]{euscript}
\usepackage{amssymb}
\usepackage{amscd}
\usepackage{epsfig}
\newcommand\ring[1]{\mathaccent23{#1}}
\def\Vr{\ring{V}} 

\newcommand{\pf}{\noindent{\em Proof: }}
\newcommand{\epf}{\hfill\hbox{\rule{3pt}{6pt}}\\}

\newcommand{\dd}{\deg_{\T}}

\newcommand{\T}{\mathcal T}

\smartqed

\begin{document}
\title{Minimum triplet covers of binary phylogenetic $X$--trees}

\author{K.~T.~Huber \and V. Moulton \and M. Steel}
\institute{KTH \at  School of Computing Sciences,
University of East Anglia \\
Norwich, UK\\
\email{K.Huber@uea.ac.uk} 
\and
VM \at  School of Computing Sciences\\
University of East Anglia\\ Norwich, UK\\
\email{V.Moulton@uea.ac.uk} 
\and
 MS \at Biomathematics Research Centre\\
University of Canterbury\\
Christchurch, NZ\\
\email{mike.steel@canterbury.ac.nz}
}
\date{Received: date/ Accepted: date}

\maketitle

\begin{abstract}
Trees with labelled leaves and with all other vertices of degree three play an important role in systematic biology and other areas of classification.   A classical combinatorial result ensures that such trees can be uniquely reconstructed from the distances between
the leaves (when the edges are given any strictly positive lengths). Moreover,  a linear number of these pairwise distance values
suffices to determine both the tree and its edge lengths.  A natural set of pairs of leaves is provided by any `triplet cover' of the tree (based on the fact that each non-leaf vertex is the median vertex of three leaves). In this paper we describe a number of new results 
concerning triplet covers of minimum size. In particular, we characterize such covers in terms of an associated graph being a 2-tree. Also, we show that minimum triplet covers are `shellable' and thereby provide a set of pairs for which the inter-leaf distance values will uniquely determine the underlying tree and its associated branch lengths. 

\end{abstract}


\section{Introduction}

Trees play a central role in systematic biology, and other areas of classification, such as linguistics.  It is often assumed that such a tree $T$ has a labelled leaf set $X$, that all vertices have degree 1 or at least three, and that there is an assignment of a positive real-valued length to each edge of $T$.  

A classical and important result from the 1960s and 1970s asserts that any such tree $T$ with edge lengths is uniquely determined  from the induced leaf-to-leaf distances between each pair of elements of $X$.  This result is the basis of widely-used methods for inferring trees from distance data, such as the popular `Neighbor-Joining' algorithm \cite{sai}. 
Moreover, when $T$ is binary (each non-leaf vertex has degree 3) then we do not require distance values for all of the $\binom{n}{2}$ pairs from $X$ (where $n=|X|$),  since just $2n-3$ carefully selected pairs of leaves suffice to determine $T$ and its edge lengths (see \cite{GLM04};  more recent  results  appear in  \cite{DHS11},  motivated by  the irregular distribution of genes across species in biological data).   

This value of $2n- 3$ cannot be made any smaller, since a binary unrooted tree with $n$ leaves has $2n-3$ edges, and the inter-leaf distances are linear combinations of the corresponding $2n-3$ edge lengths (so, by linear algebra, these values cannot be uniquely determined by fewer than $2n-3$ equations).

There is a particularly natural way to select a subset of $\binom{X}{2}$ for $T$ when $T$ is binary.  Since each non-leaf 
vertex  is incident with three subtrees of $T$, let us  (i) select a leaf from each subtree,  (ii) consider the three pairs of leaves 
we can form from this triple, and  then (iii) take the union of these sets of pairs  over all non-leaf vertices of $T$. This process produces a `triplet cover' of $T$  (defined more precisely below).   

A triplet cover need not be of this minimum size (i.e. of size $2n-3$) but in this paper we characterize when it is. Also, we show that in that case the resulting triplet cover is `shellable' which implies that the inter-leaf distances defined on these pairs uniquely determine the tree and its edge lengths.  These, and other results obtained along the way complement recent work into 
phylogenetic `lasso' sets \cite{DHS11}, \cite{HS14}, as well as a Hall-type characterization of the median function on trees in \cite{DS09}.

We begin with some definitions.

\subsection{Definitions}

Let $X$ be a finite set with $|X| \ge 3$.  We denote elements in 
$X \choose 2$ and $X \choose 3$ also
by $ab$ and $abc$, respectively, where $a,b,c \in X$ are distinct. 
We refer to the elements in ${X\choose 3}$
as {\em triples}. 

A {\em (binary) phylogenetic $X$-tree} is an unrooted tree $T=(V,E)$ which has leaf set $X$, and for which 
each non-leaf vertex is unlabelled and of degree three. 
We let $B(X)$ denote the set of binary phylogenetic $X$--trees
(two such trees are regarded as equivalent if there is a graph isomorphism 
between them that maps leaf $x$ in one tree to leaf $x$
in the other tree,  for all $x \in X$).
 In evolutionary biology, the 
set $X$ usually corresponds to 
some collection of species or taxa. 

Note that a phylogenetic $X$-tree $T$ must contain at 
least one {\em cherry} $\{a,b\}$, that is, $a$ and $b$ are adjacent
with the same interior vertex of $T$. Moreover, if $|X|>3$ then each tree 
$T \in B(X)$ has at least two cherries that are vertex disjoint from each other; if $T$ has exactly two cherries
we say it is a {\em caterpillar} tree (every tree in $B(X)$ is a caterpillar when
$|X|=4$ or $|X|=5$).  When $|X|=4$, we say that $T \in B(X)$ is a
{\em quartet}, and if the two cherries of  this tree are (say) $\{a,b\}$ and $\{c,d\}$ then 
we denote $T$ by $ab|cd$.

We  let $\Vr=\Vr(T)\subseteq V$ denote the set of $|X|-2$ interior vertices of $T$.
Given $x \in X$ where $|X|\geq 4$, we let $T-x$ denote the 
phylogenetic $(X-\{x\})$-tree which is obtained by 
removing the leaf $x$ (and its incident edge) from $T$ and suppressing the resulting
degree 2 vertex.

Suppose that $\T$ is a subset of $X \choose 2$, and $T=(V,E) \in B(X)$. 
We say that a triple in ${X \choose 3}$ {\em supports} a vertex 
$v \in \Vr$ in $T$ (relative to $\T$)
if we can select leaves $a,b,c \in X$, one from each connected
component of the graph obtained by removing $v$ and its incident edges from $T$, 
such that $ab, ac, bc \in \T$. 
We call a subset $\T \subseteq {X \choose 2}$ a 
{\em triplet cover} for $T$ if for each vertex $v\in \Vr$
there is some triple in ${X \choose 3}$ that supports $v$ (relative to $\T$). 
Note that $X=\bigcup_{A \in\T} A$ holds in this case.  
Given a non-empty subset $\T \subseteq {X \choose 2}$, we define
the {\em cover graph $\Gamma(\T)=(X,\T)$ (of $\T$)} to be the graph 
with vertex set $X$ and
edge set $\T$.  

\begin{figure}[htb]
\centering
\includegraphics[scale=1.2]{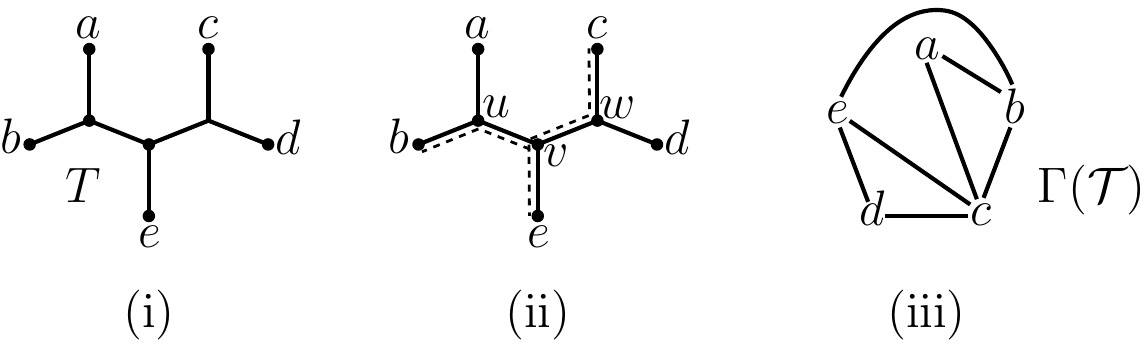}
\caption{(i) A tree $T \in B(X)$  for $X=\{a,b,c,d,e\}$; (ii) vertex $v$ is supported by the triple $bce$ (the dashed lines show the edge-disjoint paths from $v$ to these three leaves); (iii) the cover 
graph $\Gamma(\T)$ corresponding to the triplet cover $\T$ obtained by taking all pairs from  the triple $bce$ that supports $v$ and from 
the triples $abc$ and $cde$ that support vertices $u$ and $w$, respectively. This triplet cover is minimal, and since its size is 7 ($=2n-3$ for $n=|X|$) it is also a minimum triplet cover for the tree (by Proposition~\ref{prolower}). }
\label{fig1}
\end{figure}

We illustrate these concepts in Fig.~\ref{fig1}. For  the binary phylogenetic 
$X$-tree in Fig.~\ref{fig1}(i) (with $X=\{a, \ldots, e\}$) the vertex $v$ (in Fig.~\ref{fig1}(ii))  is supported 
by the triple $bce$ (there 
are three other triples that support $v$).  If $u$ is supported by, say,  
$abc$ and $w$ by $cde$ then we obtain the 
triplet cover $$\T = \binom{\{b,c,e\}}{2}  \cup \binom{\{a,b,c\}}{2} 
\cup \binom{\{c,d,e\}}{2} = \{ab, ac, bc, cd, ce, de, be\}.$$
The corresponding cover graph $\Gamma(\T)$  is shown in Fig~\ref{fig1}(iii).

Given a tree $T\in B(X)$, a  triplet cover $\T$ for $T$
is called
\begin{itemize}
\item {\em minimal} if $\T -\{ab\}$ is not a triplet 
cover for $T$, for any $ab \in \T$; 
\item {\em minimum} if
$|\T| \le |\T'|$ for every triplet cover $\T'$ for $T$.  
\end{itemize}
 These two concepts are different;  there exist 
minimal triplet covers that are not minimum (we describe an example in the final section).

Note that it can be shown that any minimum triplet cover on $X$ must have 
cardinality $2|X|-3$ (by applying Theorem~1 and Proposition~1 of  \cite{DHS11}).
Moreover, there are various ways to construct 
triplet covers that are minimum (for example, 
`pointed covers'  \cite[Theorem 7]{DHS11} and `stable triplet covers' 
\cite[Theorem 1]{HS14}).

\subsection{Outline of main results}

In this paper, we prove a structural result concerning minimum triplet covers.
Namely, we prove that a set $\T \subseteq \binom{X}{2}$ is a minimum triplet cover for a tree $T \in B(X)$ 
if and only if the associated cover graph $\Gamma(\T) = (X,\T)$ is a 2-tree (see
Theorem~\ref{cpro} and
Section~\ref{2tree} for the definition of a 2-tree). 

Using the concepts that we develop to prove this result, we also 
give an independent proof (that does not require the notion of phylogenetic `lassos' from \cite{DHS11})
that any minimum triplet cover on $X$ must have cardinality $2|X|-3$ (Proposition~\ref{prolower}). 
As a corollary of our structural result, we also 
show that if $\T$ is a minimum triplet cover for $T$ then it is shellable
for $T$ (Proposition~\ref{shellproof}).

This corollary has two important implications.
First it implies (from results in \cite{DHS11}) 
that if  $\T$ is a minimum triplet cover for $T$, then $T$ (together with its edge lengths) 
can be uniquely reconstructed from the 
tree metric restricted to the pairs in $\T$. Note that this 
can also be deduced from results in \cite{LM98} that relate 
2-trees and tree metrics (see also \cite{GLM04}).

Second, the corollary  gives an independent proof of 
\cite[Theorem 7]{DHS11} and \cite[Theorem 1]{HS14}
which state that pointed triplet covers and stable triplet 
covers are shellable, respectively.

\section{The support graph}

In this section we introduce a graph that can be associated to a 
triplet cover of a tree. Properties of this graph
will be used to help prove our results later on. 
We begin with some further definitions.

Suppose for the following that $T=(V,E) \in B(X)$.  Given a 
subset $\T \subseteq {X \choose 2}$ 
and $v \in \Vr$,  we let $S_v(\T)$ be the subset of $X \choose 3$ which 
contains precisely those triples in ${X\choose 3}$ that 
support $v$ (relative to $\T$).
We call $S_v(\T)$ the {\em support of $v$ (relative to $\T$)}.
In addition, suppose that $a,b,c\in V$ are
pairwise distinct. Then we call the unique vertex of $T$ that 
simultaneously lies on the shortest path from $a$ to $b$, from 
$b$ to $c$, and from $a$ to $c$ the {\em median} of $a$, $b$, and $c$,
denoted by ${\rm med}_T(a,b,c)$.
The following observation linking medians with supports 
will be useful.

\begin{lemma}
\label{lem1}
Let $T=(V,E)\in B(X)$ and $\T \subseteq {X \choose 2}$. 
If $abc \in S_v(\T)$, $v \in \Vr$, 
then $v = {\rm med}_T(a,b,c)$. Moreover,  $\T$ is a triplet 
cover of $T$ if and only if $|S_v(\T)| \ge 1$ for all $v \in \Vr$.
\end{lemma}

\begin{figure}[htb]
\centering
\includegraphics[scale=1.0]{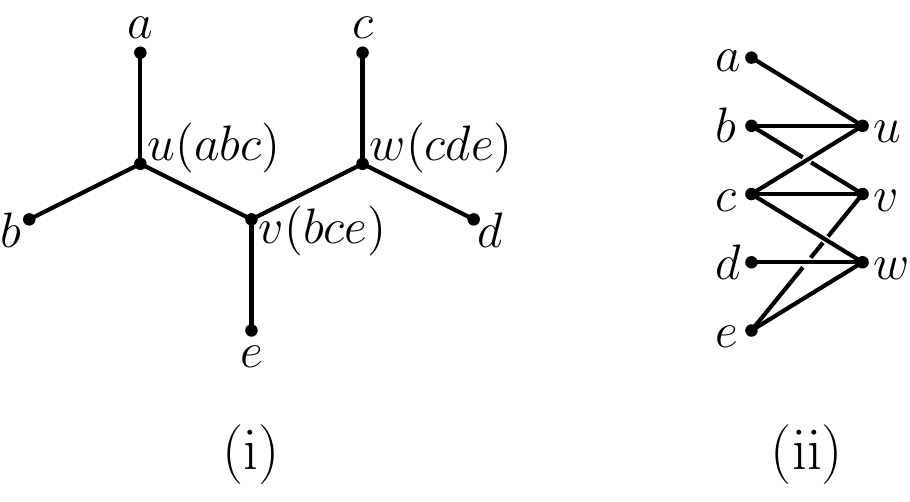}
\caption{For triplet cover $\T$ for the example from Fig.~\ref{fig1} (reproduced in (i), with the 
triple supporting an interior vertex shown in parentheses), the corresponding 
support graph $G(\T)$ is shown in (ii).}
\label{fig_support}
\end{figure}

Now, given a non-empty subset $\T \subseteq {X \choose 2}$ 
and some $x \in X$,
we put $$\T^{-x} = \T - \{ xa \,: \, a \in X-\{x\} \mbox{ and } xa \in \T\}.$$
Put differently, 
$\T^{-x}$ is the subset of $\T$ obtained by removing from $\T$ 
precisely those elements in $\T$ which contain $x$.
We also define 
a bipartite graph $G(\T) = (X \amalg \Vr, E(\T))$, with edge
$\{x,v\} \in E(\T)$, $x \in X$, $v \in \Vr$, if 
$x \in A$ for all $A \in S_v(\T)$. We call $G(\T)$ the 
{\em support graph} associated to $\T$.
For any vertex $p$ of   $G(\T)$, we let 
$\dd(p)=\deg_{G(\T)}(p)$ denote the degree of $p$ in  $G(\T)$.
In Fig.~\ref{fig_support}(ii) we illustrate the support graph 
for the triplet cover $\T$ given in Fig.~\ref{fig1}.

We now list some properties of $G(\T)$.

\begin{proposition}\label{properties}
\mbox{}
Suppose that  $\T$ and $\T'$ are triplet covers of a tree
$T=(V,E) \in B(X)$, and that $x \in X$.

\begin{itemize}

\item[(P1)] If $v \in \Vr$, then $0 \le \dd(v) \le 3$, and $1 \le \dd(x) \le |X|-2$.

\item[(P2)] If 
$\T' \subseteq \T$, then  $E(\T) \subseteq E(\T')$. In particular, 
if there exists some $x \in X$ with $\deg_{\T'}(x)=1$, then $\deg_{\T}(x)=1$.
 
\item[(P3)] If $\T$ is a minimal triplet cover for $T$,
then for all $ab \in \T$, there exists some $v \in \Vr$ such that
$a,v,b$ is a path in $G(\T)$. 

\item[(P4)] Suppose that $v$ is the vertex adjacent to $x$ in $T$. 
Then $\{v,x\} \in E(\T)$. Furthermore
$\dd(x) = 1$ if and only if $\{v,x\}$ is the only
edge in $G(\T)$ that contains $x$.

\item[(P5)]  $\dd(x) = 1$ if and only if $\T^{-x}$ is a triplet cover of 
$T-x$. 
 
\item[(P6)] 
If $\dd(x) = 1$, then $|\T| \ge |\T^{-x}|+2$.

\end{itemize}

\end{proposition}
\pf
\noindent{\em (P1):} The inequality 
$\dd(v) \le 3$ follows
immediately from the definition of 
the support $S_w(\T)$ of a vertex $w\in\Vr$ and the fact that
$T$ is binary. 
The inequality $1 \le \dd(x)$ follows 
since $x \in A$ for all $A \in S_u(\T)$
for the vertex $u$ that is adjacent to $x$ in $T$.  The inequality $\dd(x) \leq |X|-2$ follows from the fact that $T \in B(X)$ and so has
$|X|-2$ interior vertices.

\noindent{\em  (P2):} Suppose 
that $\{v,x\} \in E(\T)$, $x \in X, v \in \Vr$. Then
$x \in A$, for all $A \in S_v(\T)$. Since $S_v(\T')\subseteq S_v(\T)$
as $\T'\subseteq \T$ it follows that  $x \in A$ for all 
$A \in S_v(\T')$. Hence, $\{v,x\} \in E(\T')$. 
The second statement is a trivial consequence in light of the inequality $1 \le \dd(x)$ from (P1). 

\noindent{\em (P3):} Suppose for contradiction
that there exists some $ab\in\T$ such that for
all $v\in \Vr$, we have that $a,v,b$ is not a path
in $G(\T)$. Then for all $v\in \Vr$ there
must exist some $A\in S_v(\T)$ such that $ab\not\subseteq A$.
Hence, $\T'=\T-\{ab\}$ is a triplet
cover of $T$. Since $\T'\subsetneq \T$ clearly holds,
we obtain a contradiction in view of the minimality of $\T$.

\noindent{\em (P4):} 
That $\{v,x\}\in E(\T)$ holds is an immediate
consequence of the choice of $v$.
If $\dd(x) =1$, then since $x \in A$ for all $A \in S_v(\T)$, 
it follows that $\{x,v\}$ is in $E(\T)$.
The rest of the statement follows immediately. 

\noindent{\em  (P5):}
Suppose that $\T^{-x}$ is not a triplet cover of $T-x$. Then, 
by Lemma~\ref{lem1}, there exists 
an interior vertex $u$ of $T-x$ such that $S_u(\T^{-x})=\emptyset$.   
Let $u' $ be the vertex in $T$ that corresponds to $u$ in $T-x$. Then
as  $S_u(\T^{-x})=\emptyset$, it 
follows that $x \in A$ for all $A \in S_{u'}(\T)$.
Hence $\{x,u'\}\in E(\T)$ and, so, $\deg_{\T}(x)\geq 1$. 
Moreover, if $v$ is the vertex 
adjacent to $x$ in $T$, then $v \neq u'$. By (P4), it follows that
$\{x,v\}$ is also an edge in $E(\T)$.
Therefore $\deg_{\T}(x)>1$.

Conversely, suppose that $\T^{-x}$ is a triplet cover for $T-x$, 
and assume for contradiction that $\deg_{\T}(x) \ge 2$. 
Then there exist $u,v \in \Vr$ distinct  
such that $x \in A$ for all $A \in S_u(\T)$ and $x \in B$ for all
$B \in S_v(\T)$. Without loss of 
generality, we may assume
that $v$ is the vertex in $T$ that is 
adjacent to $x$. Let $u'$ be
the vertex in $T-x$ that corresponds to $u$ 
in $T$. Then $S_{u'}(\T^{-x})=\emptyset$
since $x\in A$ for all $A \in S_u(\T)$. 
Hence $\T^{-x}$ is not a triplet cover for
$T-x$, a contradiction.  

\noindent{\em (P6):}
If $v$ is the vertex in $T$ 
adjacent to $x$, then  $S_v(\T) \neq \emptyset$ by Lemma~\ref{lem1}.
Hence, there must be some $A \in S_v(\T)$ 
with $x \in A$. But then $|\T - \T^{-x}|\ge 2$. 
\epf

We now show that any minimal triplet  cover of a tree in $B(X)$ has a size that grows linear with $|X|$.

\begin{corollary}\label{M0}
Suppose that $\T$ is a minimal triplet cover of
some  $T\in B(X)$. 
Then 
$$
|\T| \le 3(|X|-2).
$$
\end{corollary}
\pf
Put $T=(V,E)$.
First we observe that if $B=(X \amalg \Vr,E')$ is a bipartite graph
in which every vertex in $\Vr$ has degree at most 3, then 
the number of length 2 paths in $B$ of the form
$x, v, y$ with $x,y \in X$ and $v \in \Vr$  is equal to 
$$
\sum_{v \in \Vr} | \{x, v, y \,:\, x,y \in X \mbox{ and } 
x, v, y \mbox{ a path in } B \}|.  
$$
Now, by 
(P3), $|\T|$ is 
less than or equal to the number of length 2 paths in 
$G(\T)$ of the form $x, v, y$ with $x,y \in X$ and $v \in \Vr$.
Since $|\Vr|=|X|-2$, and each term in the above sum is at most 3 the corollary follows.
\epf

\section{Multiplicities}

In this section we derive some bounds for degrees of vertices in the cover
graph of a triplet cover. 
Suppose that $\T$ is a triplet cover of $T \in B(X)$. 
For $x \in X$ we define the {\em multiplicity $\mu(x)=\mu_{\T}(x)$
of $x$ (relative to $\T$)} to be the number of elements in $\T$
that contain $x$ (or in other words, the degree of the vertex
$x$ in the cover graph $\Gamma(\T)$).
The {\em multiplicity of $\T$} is $\mu(\T) = \min\{\mu_{\T}(x): x\in X\}$.

The following observation relating multiplicities with 
degrees will be useful later.

\begin{lemma}\label{M4}
Suppose that $\T$ is a triplet cover for some tree $T\in B(X)$
and $x \in X$. 
If $\mu(x)=2$, then $\dd(x)=1$.
\end{lemma}
\pf
If $\mu(x)=2$, then $x$ can be contained
in at most one element of $\bigcup_{v \in \Vr} S_v(\T)$. 
But $x$ must be contained
in every element of $S_u(\T)$ for $u$ the vertex in $\Vr$ that is adjacent to $x$ in $T$.
Hence $|S_u(\T)|=1$, and the only edge contained in the support graph $G(\T)$ that
contains $x$ (which must exist by (P1)) is $\{x,u\}$. In particular,  $\deg_{\T}(x)=1$.
\epf

We now derive some bounds for multiplicities of minimal and minimum triplet covers.

\begin{proposition}\label{multbounds}
Suppose that $T\in B(X)$.
\begin{itemize}

\item[(M1)] If $\T$ is a minimal triplet cover for $T$, then $2 \le \mu(\T) \le 5$.

\item[(M2)] If $\T$ is a minimum triplet cover for $T$, then $2 \le \mu(\T) \le 3.$

\end{itemize}

\end{proposition}

\pf \noindent{\em (M1):}
Suppose that $x \in X$. Let $v$ be the vertex in $T$ 
adjacent to $x$ in $T$. Then, as $\T$ is a triplet cover for $T$,
by Lemma~\ref{lem1} there must exist some $axy \in S_v(\T)$
where $a,y \in X-\{x\}$ are distinct. 
Therefore  $2 \le \mu(x)$ for all $x \in X$ and
so $2 \le \mu(\T)$. 

To see that the remaining inequality holds, we show that 
there is some element of $X$ that is
contained in at most 5 elements of $\T$.
We use a simple counting argument based on pairs
$(x,c)$ where $x \in X$ is an
element in some $c \in \T$. By Corollary~\ref{M0},
$|\T|\leq 3(|X|-2)$ as $\T$ is minimal. Since
each element of $\T$
contains 2 elements of $X$, the size of the 
set $R$ of pairs
$(x,c)$ is at most $6(|X|-2)$. On the other hand
$\sum_{x \in X} \mu(x) = |R|$.
Hence, since $|X|\geq 3$, there must exist some $x \in X$ with 
$\mu(x) \leq 5$. 

\noindent{\em (M2):}
We again count pairs $(x, c)$ where $x \in c$ and $c$ is an element in $\T$. 
This is $2 |\T | = 2(2|X|-3)$ and also equal to
$\sum_{x \in X} \mu(x)$.
Since $2(2|X|-3)<4|X|$ and $|X|\geq 3$, 
there is some $x \in X$ with $\mu(x) \leq 3$.
That $\mu(\T)\geq 2$ holds follows from (M1).
\epf

\section{A lower bound}
In this section, we show that a minimum triplet cover of 
a tree $T \in B(X)$ has size $2|X|-3$. As mentioned
in the introduction, this result can also be derived 
by applying Theorem~1 and Proposition~1 of  \cite{DHS11}. However, it is of 
interest to have a direct proof that is independent of results concerning tree metrics.

\begin{proposition}
\label{prolower}
\mbox{}
Suppose that $\T$ is a triplet cover for some $T \in B(X)$. Then 
$|\T| \ge 2|X|-3$. Moreover this bound is tight.
\end{proposition}

\pf We use induction on $n =|X|$.
The result clearly holds for $n=3$. So, suppose that the result holds
for all triplet covers of trees in $B(X)$
with $3 \le |X| \le n-1$.

Suppose that $\T$ is a triplet cover for a tree in $B(X)$ with $|X|=n$. 
If there exists some
$a \in X$ such that $\dd(a)=1$, then by (P5) $\T^{-a}$
is a triplet cover for $T-a$. Hence, by (P6) and
induction, $|\T| \ge |\T^{-a}| + 2 \ge 2n-3$.

So, suppose that  $\dd(x)\ge2$ for all $x \in X$. 
Note that there
must exist some $a \in X$ with $\dd(a)=2$ (otherwise, 
$\dd(x)\ge3$ for all $x \in X$ implies 
that there is a vertex $v \in \Vr$ with $\dd(v)\geq 4$, which contradicts (P1)).
Suppose that $v,u \in \Vr$ are distinct with $\{a,v\}, \{a,u\}$
in $E(\T)$. Then there exist distinct elements $b,c,x,y \in X-\{a\}$  
with $\{b,x\} \neq \{c,y\}$ such that
$abx \in S_v(\T)$ and $acy \in S_u(\T)$. Put 
$C:=\{b,x\}\cap \{c,y\}$. Then since  $\{b,x\} \neq \{c,y\}$ it follows that $|C|<2$ and so  we 
consider the two possible cases ($|C|=1$ and $|C|=0$).

\noindent{\em Case 1:} $|C|=1$. Without loss of 
generality we may assume
$x=c$ and $y \neq b$. 
Then it is straight-forward to see that without
loss of generality,  $v$
is adjacent to $a$ in $T$, $u$ lies on the path in $T$
between $v$ and $c$, and
$T$ restricted to the
set $\{a,b,c,y\}$ is the quartet $ab|cy$.
Note that $by \not\in \T$ since otherwise $bcy \in S_u(\T)$
which contradicts $\{a,u\} \in E(\T)$. 

Consider the triplet cover $\T' = \T \cup \{by\}$ of $T$.
Then $acy,bcy \in S_u(\T')$. Hence, since $E(\T')\subseteq E(\T)$ 
by (P2), $\deg_{\T'}(a)=1$.
Therefore, by (P5), $\T'^{-a}$ is a triplet
cover of $T-a$. But the elements $ab, ac, ay$ of $\T$ 
are not contained in $\T'^{-a}$ and, so,
$$
|\T'^{-a}|+3 \le |\T'| = |\T|+1. 
$$
The fact that $ |\T| \ge 2|X|-3$ holds now follows immediately 
by induction.

\noindent{\em Case 2:} $|C|=0$. Then 
$x \neq c$ and $y \neq b$. Without loss of generality, 
we can assume that $v$ is adjacent to $a$ in $T$, and
that $T$ restricted to the set $\{a,b,c,y,x\}$ is a caterpillar tree with cherry $\{a,x\}$.
We consider the case where $\{y,c\}$ is also a cherry in this caterpillar tree and $u$ 
is adjacent to both $y$ and $c$ in $T$. 
The argument for the remaining case
(where $\{b,y\}$ or $\{b,c\}$ is also a cherry) is similar.

First note that if $bc \in \T$, then $by \not\in \T$, since otherwise
$byc \in S_u(\T)$ which would contradict $\{a,u\} \in E(\T)$. 
Similarly if $cx \in \T$, then $yx \not\in \T$. Hence, by symmetry,
we can assume that $\T$ does not contain at least one 
element from the set $\{bc, by\}$ and at least one element from the set $\{cx,yx\}$.
Now, let $P$ be a subset of $\{bc,by,cx,yx\}-\T$ of minimum size such that 
$\T \cup P$ contains precisely one of the sets $\{bc, by\}$ or $\{cx,yx\}$, noting that $|P|\le 2$.
Consider the triplet cover $\T' = \T \cup P$ of $T$.
Then it is easily seen that $\deg_{\T'}(a)=1$, and so by (P5) $\T'^{-a}$ is a triplet
cover of $T-a$. But the elements $ab, ac, ax, ay$ of $\T$ 
are not contained in $\T'^{-a}$ and so
$$
|\T'^{-a}|+4 \le |\T'| = |\T|+|P| \le |\T|+2.
$$
The fact that $|\T| \ge 2|X|-3$ holds now follows by induction.

The fact that the bound is tight follows since 
for every $T \in B(X)$ there exists 
some triplet cover of $T$ with cardinality $2|X|-3$ 
(e.g. a pointed cover \cite{DHS11}).
\epf

\section{A characterization of minimum triplet covers} \label{2tree}

In this section, we prove our main result, namely a characterization
of minimum triplet covers in terms of 
the structure of their cover graphs. First, we recall that 
a graph $H=(V,E)$ is called a {\em 2-tree} if 
there exists an ordering $v_1,v_2,\dots,v_m$ of $V$ such
that $\{v_1,v_2\} \in E$ and, for $i=3,\dots,m$, the 
vertex $v_i$ has degree 2 and belongs to a 
unique triangle in the subgraph induced by $H$ on
the set $\{v_1,v_2,\dots,v_i\}$ \cite[p.235]{GLM04}.
It is easily seen that a 2-tree has treewidth at most 2, and 
conversely, every graph of treewidth at most 2 is a subgraph of a 2-tree. 


\begin{theorem}
\label{cpro}
\mbox{}
Suppose that $\T$ is a triplet cover for a tree $T \in B(X)$.
Then $\T$ is minimum triplet cover if and only if 
$\Gamma(\T)$ is a 2-tree.
\end{theorem}
\pf Put $T=(V,E)$.
Suppose that $\Gamma(\T)$ is a 2-tree. Then since
2-trees on $n$ vertices have $2n-3$ edges \cite[p.227]{LM98}
and $|X|=n$, we have 
$\T = 2|X|-3$. So $\T$ is a minimum triplet cover for 
$T$.

Conversely, suppose that $\T$ is a minimum triplet cover
for some tree $T \in B(X)$.
We shall prove that $\Gamma(\T)$ is a 2-tree by induction on $n=|X|$.
If $|X|=3,4$ it is clearly true. Suppose the statement holds
for all $X$ with $3 \le |X| \le n-1$.
 
Let $\T$ be a minimum triplet cover for $T$ on $X$ with $n=|X|$. 
Note that, by (M2), $\mu(\T)$ 
equals 2 or 3. Also, note that $\T$ must be a minimal 
triplet cover for $T$.

\begin{figure}[htb]
\centering
\includegraphics[scale=1.0]{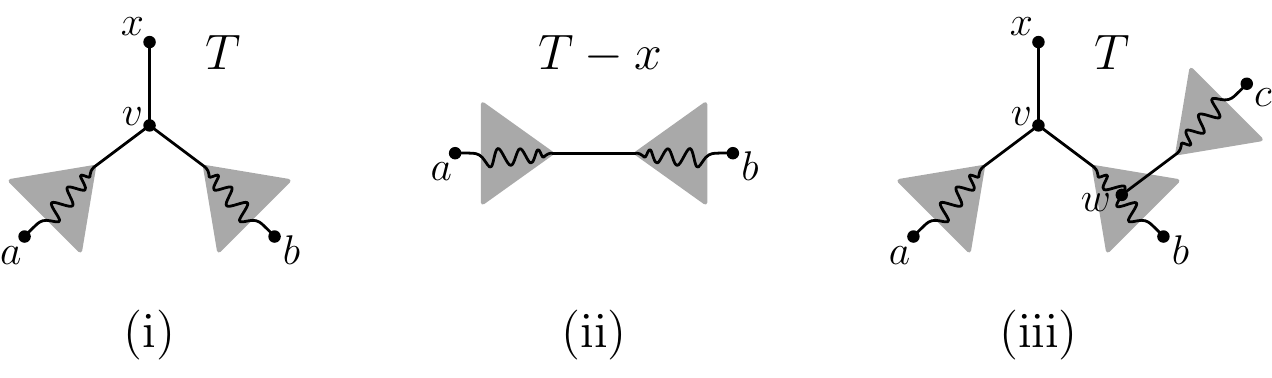}
\caption{ Figures for the proof of Theorem~\ref{cpro}. (i) Leaf $x$  and the other two leaves that form the triple in  $S_v(\T)$; (ii) the tree $T-x$ obtained from $T$ by restricting this tree to $X-\{x\}$; (iii) The labelling of additional vertices in the case where $\mu(\T)=3$. Squiggly lines denote paths in $T$.}
\label{fig_extra}
\end{figure}

Suppose that $\mu(\T)=2$. Let $x \in X$ be such that $\mu(x)=2$. 
Then there exist
$a,b \in X-\{x\}$ with $xa, xb \in \T$. Consider the vertex 
$v\in V(T)$ adjacent 
to $x$ in $T$ (as shown in Fig.~\ref{fig_extra}(i)).  Then as $\T$ is a 
triplet cover, and $xa, xb$
are the only elements in $\T$ containing $x$, it follows that 
$S_v(\T) = \{xab\}$. 
Hence, $ab \in \T$. 
It follows that $\T':= \T-\{xa, xb\}$ is a triplet cover for $T-x$ (see Fig.~\ref{fig_extra}(ii))  and since $|\T|= 2|X|-3$, it follows that
$|\T'| = 2(|X|-1)-3$ and so $\T'$ is a minimum triplet cover for $T-x$. Since $T-x$ has one fewer leaf than $T$, we can apply the induction hypothesis and conclude that
$\Gamma(\T')$ is a 2-tree. Then, since $\Gamma(\T)$ is obtained from $\Gamma(\T')$ by
attaching $x$ to the endpoints of the edge $\{a,b\}$ in $\Gamma(\T')$, it follows that 
 $\Gamma(\T)$ is also 2-tree.

Now suppose  that $\mu(\T)=3$. We shall show that this is 
not possible, from which the theorem follows. Let $x \in X$ 
be such that $\mu(x)=3$ and let $v \in \Vr$ denote 
the vertex adjacent to $x$ in $T$. 
Then since $\T$ is a minimal triplet cover for $T$
there must exist $a,b \in X-\{x\}$ distinct such
that $xab\in S_v(\T)$. Moreover, as $\mu(x)=3$
there must exist some $c \in X-\{x,a,b\}$ with $xc \in \T$.
Since we also have $xa,xb\in \T$, and since $\T$ is a minimum triplet cover, it follows that $bc\in \T$.

Without loss of generality, assume $T$ restricted to $x,a,b,c$
is the quartet $xa|bc$ (notice that we have symmetry involving $a$ and $b$, and the quartet cannot be $xc|ab$ because 
of the assumption that $xab \in S_v(\T)$ where $v$ is the vertex adjacent to $x$ in $T$), as shown in Fig.~\ref{fig_extra}(iii). 
Let $w \in \Vr$ be such that $w = {\rm med}(x,b,c)$.

We claim that $ac \not\in \T$. 
Assume for contradiction
that $ac\in\T$.
Since $\T$ is minimal and $xc\in\T$, there exists some vertex 
$u\in\Vr$ and some $A\in S_u(\T)$ such that $xc \subset A$.
Note that as $\mu(x)=3$, we must have $u \in \{v,w\}$.
If $u=v$ then $\T-\{xc\}$  is a smaller minimum triplet cover for $T$ (since $v$ is still supported by $abx$), and this contradicts the minimality of $\T$. 
Thus we may assume that $u=w$, in which case there is a set $A \in S_w(\T)$ with $xc \subset A$.  Since $\mu(x)=3$ and we already have $ax, bx, cx \in \T$ it follows that
$A=xbc \in S_w(\T)$ which implies that $bc \in \T$.  However, as we already have $ab \in \T$,  the additional assumption that
$ac \in \T$  means that $\T-\{xc\}$ contains $ab, ac, bc$ which provides an alternative set, namely $abc$ in $S_w(T)$, in which case $\T-\{xc\}$ remains a triplet cover for $T$. But again this contradicts the minimality of $\T$.
Thus, $ac \not\in \T$, as claimed. 

Therefore, in summary, $xa,xb,xc,ab,bc \in \T$ and $ac \not\in \T$.
We claim next that $\T' = \T - \{xb\} \cup  \{ac\}$
is a triplet cover for $T$. Indeed, 
 if $xb$ is contained in some
element in $S_u(\T)$ for some $u \in \Vr$, then since 
$\mu(x)=3$ we must have $u \in \{v,w\}$. Since
$acx \in S_v(\T')$ and $abc \in S_w(\T')$ it follows that 
$\T'$ must be a triplet cover for $T$, as claimed.

To complete the proof, note that
since $\mu_{\T'}(x)=2$, Lemma~\ref{M4} implies  $\deg_{\T'}(x)=1$. 
Hence, by (P5), $\T'^{-x}=\T' - \{xa,xc\}$ is a triplet 
cover of $T-x$. 
Since $T-x$ has one fewer leaf than $T$ we can apply the induction hypothesis and conclude that the graph $\Gamma(\T'^{-x}) = (X-\{x\},\T'^{-x})$ 
is a 2-tree. Since any 2-tree has at least two vertices
with degree 2 \cite[p.227]{LM98}, it follows
that in $\Gamma(\T'^{-x})$
at least
one of the two vertices $a$ or $c$ has degree 2
(since there cannot be a vertex $y \in X-\{x,a,b,c\}$
such that the degree of $y$ in $\Gamma(\T'^{-x})$ is equal to 2 as, by assumption,
$\mu(\T)=3$). But if, without loss
of generality, the degree of $a$ in $\Gamma(\T'^{-x})$ is equal to 2, 
then $\mu_{\T}(a)=2$ must hold too
which contradicts $\mu(\T)=3$. This completes the 
proof.
\epf

The next result follows immediately from the last theorem and the fact
that any 2-tree has at least two vertices
with degree 2 (see e.g. \cite[p.227]{LM98}). It 
improves on the bound given in
Proposition~\ref{multbounds} (M2).

\begin{corollary} \label{C2}
If $\T$ is a minimum triplet cover for some tree
$T \in B(X)$ then $\mu(\T)=2$.
\end{corollary}

Note that a 2-tree is a 2d-tree, but not necessarily conversely.   
\cite[Proposition 3.4]{GLM04} (a graph $G = (V, E)$ is called a {\em 2d-tree}  if there exists
an ordering $x_1, x_2, \ldots, x_n$ of $V$ such that $\{x_1, x_2\} \in E$ and, for $i=2, \ldots, n$
the  vertex $x_i$ 
has degree 2 in the subgraph of $G$ induced by $\{x_1, x_2, \ldots, x_i\}$).
So Theorem~\ref{cpro}  can be used to strengthen
Theorem 1 of  \cite{HS14}.

\section{Shellings}
\label{shell}

Given a triplet cover $\T$ of a 
tree $T \in B(X)$, we say that $\T$ is 
{\em $T$-shellable} if there exists an ordering of the elements in 
${X \choose 2} - \T$, say $a_1b_1, a_2b_2, \dots, a_mb_m$
such that for every $1 \le i \le m$, there exists
a pair $x_i,y_i$ of distinct elements in $X - \{a_i,b_i\}$
such that the restriction of $T$ to the 
set $Y_i = \{a_i,b_i,x_i,y_i\}$ is the
quartet $x_ia_i | y_ib_i$, and all elements in $Y_i \choose 2$
except $a_ib_i$  are
contained in $\T_i = \T \cup \{a_jb_j \,: 1 \le j \le i-1\}$. If $T$ is
clear from the context then we sometimes
just say that $\T$ is {\em shellable}, and 
we refer to the ordering of ${X \choose 2} - \T$ as a {\em shellable ordering}.

Although this combinatorial definition of shellability  seems somewhat involved, its motivation rests on it being a sufficient condition for recursively determining the distances between all pairs of leaves (when the edges of $T$ are assigned arbitrary positive edge lengths)  starting with just the distance values for the pairs in the triplet cover. In other words,  if a  triplet cover $\T$ of a tree $T \in B(X)$ is shellable then the pairs of elements from $X$ that are not already present in $\T$ can be ordered in a sequence so 
that the distance in $T$ between the leaves in each pair is uniquely determined from the distances values on pairs that are either (i) present
as an element of $\T$ or (ii) appear earlier in the sequence.

For example, for the tree $T$ shown in Fig.~\ref{fig1}(i),  and the triplet cover $\T$ consisting of the 7 pairs of elements of $X$ that form the edges of $\Gamma(\T)$ in Fig.~\ref{fig1}(iii), there are just three pairs from $\binom{X}{2}$ that are not present
in $\T$, namely $ad, ae, bd$.   Ordering the  pairs as $a_1b_1 = ae, a_2b_2=ad, a_3b_3=bd$ provides a shellable ordering, since for $ae$ we can select $x_1y_1 = bc \in \T$ and observe that  
$x_1a_1|y_1b_1 = ba|ce$ is the quartet obtained by restricting $T$ to $\{a,b,c,e\}$, the distance between $a_1=a$ and $b_1=e$ in $T$ is determined uniquely by the five
other distances involving pairs from $\{a,b,c,e\}$, and these five pairs are present in $\T$.  Having determined the distance for $a_1b_1$ one can now use this (and the distances
for pairs in $\T$) to compute the distance value for the pair $a_2b_2$ and, subsequently, for the pair  $a_3b_3$.

We now gather together some facts concerning the shellability of 
triplet covers, including shellability of minimum triplet covers.

\begin{proposition}\label{shellproof}
\label{spro}
\mbox{}
\begin{itemize}
\item [(S1)] Suppose that $T\in B(X)$,
$x\in X$, and $\T$ is a triplet cover of $T$ 
such that $\T^{-x}$ is a triplet cover of $T-x$. 
If $\T^{-x}$ is ($T-x$)-shellable, then $\T$ is $T$-shellable.

\item[(S2)] Suppose that $\T$, $\T'$ are triplet
covers of some tree $T \in B(X)$ and
that $\T' \subseteq \T$. If $\T'$ is 
$T$-shellable, then so is $\T$.

\item[(S3)] If $\T$ is a minimum triplet cover for a tree $T \in B(X)$,   
then $\T$ is $T$-shellable.

\end{itemize}
\end{proposition}

\pf {\em (S1):} Put $T=(V,E)$. Suppose $x \in X$ such that
$\T^{-x}$ is a triplet cover of $T-x$ which is shellable. Suppose that
$v \in \Vr$ 
is the vertex in $T$ that is adjacent to $x$ in $T$. Then 
there must exist $a,b \in X-\{x\}$ distinct with $xab\in S_v(\T)$.
Let $\T(x) =  \{ de \in \T \,: \, x \in \{d,e\} \}$
and $\T^*(x) =  \{ de  \in {X \choose 2} \,: \, x \in \{d,e\}
 \mbox{ and } de \not\in \T(x) \}$, so that
$\T = \T^{-x} \amalg \T(x)$ and 
$$
{X \choose 2} - \T = \left( {X-\{x\} \choose 2} - \T^{-x} \right)
 \amalg \T^*(x).
$$
Since $\T^{-x}$ is $(T-x$)-shellable, there 
is a shellable ordering of ${X-\{x\} \choose 2} - \T^{-x}$ so that 
all of the elements in that set can be added into $\T^{-x}$ 
to obtain ${X-\{x\} \choose 2}$. 

To complete the shellable ordering it remains to add the elements of $\binom{X}{2}$ that contain $x$ to the ordering so far constructed. 
We consider two cases. First, suppose 
that neither 
$\{x,a\}$ nor $\{x,b\}$ form a cherry of $T$. 
Then for all $px \in \T^*(x)$, without loss of generality, 
the quartet induced by $T$ on $\{x,a,b,p\}$ 
is $ap|xb$. Since we have that
$xa,xb,ab$ as $xab\in S_v(\T)$ and also $ap$ and $bp$ as we
have all elements in ${X-\{x\} \choose 2}$, it follows that we can add 
in $xp$ as a next element of the shellable ordering.  
We can repeat this adding-in process
for all remaining elements in $\T^*(x)$ (in any order) to 
obtain ${X \choose 2}$. So $\T$ is $T$-shellable in this case.

Second, suppose without loss of generality that $\{x,a\}$ forms a cherry.
Then if $px \in  \T^*(x)$, then the quartet induced by $T$ 
on the set $\{x,a,b,p\}$ is $xa|bp$.
So, using similar arguments as in the
previous case, we can add in $xp$ as a next element in the shellable ordering. 
It follows that  we
can repeat this process for all remaining elements in 
$\T^*(x)$ (in any order) to 
obtain a shellable ordering of ${X \choose 2}$. So $\T$ is $T$-shellable in this case too. 

\noindent{\em (S2):} This follows immediately from the
definition of shellability. 

\noindent{\em  (S3):}
We proceed using induction on $n = |X|$. For $n=4$
the statement is clearly true. Suppose 
the statement is true up to and including $n-1 \ge 4$.  

Let $\T$ be a triplet cover for some phylogenetic $X$-tree 
with $|X|= n$.
By Corollary~\ref{C2}, $\mu(\T) =2$. Suppose that $x \in X$ 
with $\mu(x)=2$. Then, by Lemma~\ref{M4}, $\dd(x)=1$.
By (P4) it follows that $\T^{-x}$ is a triplet cover for $T-x$.
Note that  $\T^{-x}$ is minimum
since $|\T^{-x}|= |\T|-2$. 
Thus by induction $\T^{-x}$ is $(T-x$)-shellable. Therefore, $\T$
is $T$-shellable by (S1).\epf

\begin{corollary}
For any tree $T \in B(X)$, suppose that $\T$ is a minimum triplet cover 
for $T$. Consider any assignment of strictly positive lengths to the 
edges of $T$, and the resulting assignment of inter-leaf distances on 
the pairs from $\T$. This function from $\T$ to ${\mathbb R}^{>0}$ 
uniquely determines $T$ and its edge lengths, since no different tree 
$T' \in B(X)$ can induce the same inter-leaf distances on pairs from 
$\T$ under any positive weighting of the edges of $T'$.
\end{corollary}
\pf
This follows immediately from Part (S3) of Proposition~\ref{shellproof}, combined with Theorem 6 of \cite{DHS11}.
\epf

Note that there are examples of sets $\T \subseteq \binom{X}{2}$ 
having cardinality $2|X|-3$ that determine $T$ and any set of positive 
edge lengths from inter-leaf  distances, but which are not $T$-shellable 
(see Example~\ref{L7}).

\begin{example} \label{L7}
Put $X=\{a,b,c,d,e,f,g\}$ and let $T$ be the  caterpillar tree
with exactly two cherries $\{a,b\}, \{f,g\}$ 
and intermediate leaves $c,d,e$ (as shown in Fig.~\ref{fig2}(ii)).
Put $\T =\{ab,ad,bc,be,cd,cf,de,dg,ef,fg,ag\}$. Then $\T$ 
determines $T$ and any set of positive edge lengths from inter-leaf 
distances, but it is not $T$-shellable
\cite[Example 6.2]{DHS11}.
\end{example}

\section{Conclusion and open problems}

As mentioned earlier, there are examples of minimal triplet covers $\T$ 
that are not minimum. The following provides a specific example.  

\begin{figure}[htb]
\centering
\includegraphics[scale=0.9]{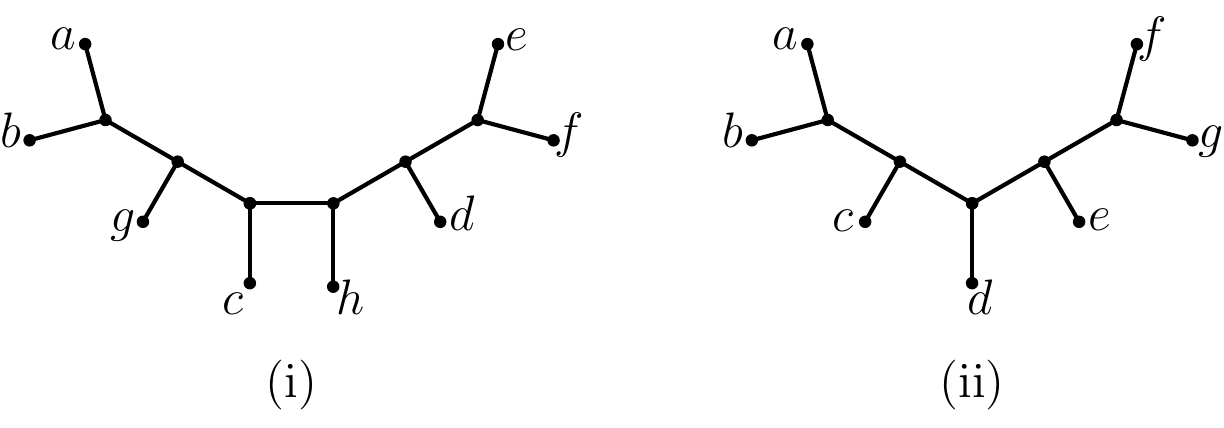}
\caption{(i) A phylogenetic $X$-tree with $X=\{a,\ldots, h\}$.
The set $\T =\{ab, ac, bc, cd, bd, ce, de, df, ef, ah, ag, fg, fh, gh\}$ is a minimal triplet cover but not a minimum one.  The set
$S_v(\T)$ associated with each interior vertex $v$ of $T$ generates the following sequence (from left-most to right-most interior vertex):
$abc, fgh, cbd, fgh, dec, edf$. 
 (ii) A  phylogenetic $X$-tree with $X=\{a\ldots, g\}$ for which the set 
$\T =\{ab,ad,bc,be,cd,cf,de,dg,ef,fg,ag\}$ determines $T$ along with an assignment of positive edge lengths from the induced inter-leaf distances, yet $\T$ is not shellable. }
\label{fig2}
\end{figure}

\begin{example} \label{minimal}
Let $X=\{a,b,c,d,e,f,g,h\}$ and
$T$ be the phylogenetic $X$-tree having  
 cherries $\{a,b\}, \{e,f\}$ and
leaves, starting with cherry $\{a,b\}$, labeled in the order 
$g,c,h,d$ (see Fig.~\ref{fig2}(i)).
Let $$\T =\{ab, ac, bc, cd, bd, ce, de, df, ef, ah, ag, fg, fh, gh\}.$$ 
Then $\T$ is a minimal triplet cover for $T$. Since $|\T|=14\not=2|X|-3$
it follows that $\T$ is not minimum.
\end{example}

\noindent An interesting problem would be to investigate the structure of the cover graph
for minimal triplet covers.

\bigskip

\noindent Our results also suggest further questions for future work.
\begin{itemize}
\item[(i)] There are formulae for counting the number of labeled 2-trees
\cite{M69}. Is there a formula for counting the number of minimum 
triplet covers for a given phylogenetic $X$-tree? 

\item[(ii)] We have shown that minimum triplet covers are shellable. It would be interesting to see how far this result extends.   For example, is {\em every} triplet cover shellable?  Understanding the
structure of minimal triplet covers  might help to shed light on this question. 
\end{itemize}

\bigskip

\noindent{\bf Acknowledgment} We thank the two anonymous reviewers for helpful comments,  particularly Reviewer 1 for numerous helpful suggestions.
 MS thanks the Allan Wilson Centre for helping fund this work. 
VM and KTH thank the London Mathematical Society for
helping to fund their visit to the University of Canterbury,
Christchurch.

\end{document}